\documentclass{birkjour}
\newtheorem{thm}{Theorem}[section]
 \newtheorem{cor}[thm]{Corollary}
 \newtheorem{lem}[thm]{Lemma}
 \newtheorem{prop}[thm]{Proposition}
 \theoremstyle{definition}
 \newtheorem{defn}[thm]{Definition}
 \theoremstyle{remark}

 \numberwithin{equation}{section}
\usepackage{amsmath,amsfonts,amssymb,amsthm,newlfont}

\begin{document}
\title[Conformal Killing vector fields and Rellich type identities]
 {Conformal Killing vector fields and Rellich type identities on Riemannian manifolds, II }
\author[Yuri Bozhkov]{Yuri Bozhkov}

\address{%
Instituto de Matem\'atica, Estatistica e Computa\c c\~ao Cient\'\i fica - IMECC \\
Universidade Estadual de Campinas - UNICAMP\\
Rua S\'ergio Buarque de Holanda, 651\\
 $13083$-$859$ - Campinas - SP\\
Brasil}

\email{bozhkov@ime.unicamp.br}

\thanks{This work was completed with the partial support of FAPESP and
CNPq, Brasil (Y.B.)}
\author{Enzo Mitidieri}
\address{Dipartimento di Matematica e Informatica \\
Universit\`a degli Studi di Trieste\\
Via Valerio $12/1$, $34127$ Trieste\\
Italia}
\email{mitidier@units.it}

\begin{abstract}
We propose a general Noetherian approach to Rellich integral identities.
 Using this method we obtain a higher order Rellich type identity involving the
 polyharmonic operator on Riemannian manifolds admitting homothetic transformations. Then
 we prove a biharmonic Rellich identity in a more general context. We also establish a nonexistence result for 
 semilinear systems involving biharmonic operators.
\end{abstract}

\newcommand{\bb}{\begin{equation}}
\newcommand{\ds}{\displaystyle }
\newcommand{\ee}{\end{equation}}
\def\am{\varpi}
\def\om{\omega}
\def\ga{\gamma}
\def\Ga{\Gamma}
\def\sig{\sigma}
\def\Sig{\Sigma}
\def\al{\alpha}
\def\Te{\Theta}
\def\eps{\varepsilon}
\def\ups{\upsilon}
\def\var{\varphi}
\def\La{\Lambda}
\def\la{\lambda}
\def\de{\delta}
\def\Om{\Omega}
\def\te{\theta}
\def\th{\vartheta}
\def\be{\beta}
\def\De{\Delta}
\def\r{\rho}
\def\ka{\mbox{\large $\kappa$}}
\def\ze{\zeta}
\def\l1{{\lambda}_1}

\def\lm{{\Delta }^m_g}

\def\nb{\nabla}
\def\kd{\partial}

\def\R{{\rm I\!R}}
\def\nb{\nabla}
\def\F{{\cal F}}
\def\A{{\cal A}}

\maketitle

\section{Introduction}

It is well-known that integral identities play an important role
in the theory of functions and differential equations. Basically
 there are two kinds of such identities which are most frequently
 used, namely, Pohozaev type identities \cite{po01,po02} and
 Rellich type identities \cite{rel}. An important observation to be  made is
 that the Pohozaev identities are satisfied by solutions of
 Dirichlet boundary problems while the Rellich identities concern
 functions which belong to certain functions spaces without any
 reference to other relations which they may satisfy like
 differential equations or boundary conditions. For this reason, as it has been pointed out
 in \cite{ye2},  the Rellich's Identity is an `important tool for obtaining,
  among other things,
 a priori bounds of solutions for semilinear Hamiltonian elliptic systems \cite{cfm}, nonexistence results \cite{em0,em} and sharp Hardy type inequalities
 \cite{em2}.

 The main purpose of this paper is to propose an unified approach
 to both Rellich and Pohozaev type identities. We have initiated
 this research with the paper \cite{ye} in which we devised and
 developed a Noetherian approach to Pohozaev's identities
 whose essential point is that the latter can be obtained from
 the Noether's identity \cite{i1,i2} after integration and
 application of the divergence theorem taking into account the corresponding
 equations (or systems) and boundary conditions. In this procedure one chooses critical
 values of the involved parameters. (See \cite{ye1} on the last
 point.) In a subsequent work \cite{ye2} we have applied this method to
 some semilinear partial differential equations and systems on
 Riemannian manifolds. For this
 purpose we have employed conformal Killing vector fields, which
 are related to Noether symmetries of critical differential equations (see
 \cite{yi}). In fact, this generalizes the original idea of
 Pohozaev \cite{po01} who has made use of the radial vector filed $X=r\frac{\kd }{\kd r}$ on
 $\mathbb{R}^n$. We have also obtained in \cite{ye2} a Rellich
 type identity on manifolds following the argument in
 \cite{em0,em}.

 The corner stone of this work is the observation that
 the Rellich type identities for functions on Riemannian manifolds
 can be generated by integration of
 the Noether's Identity for appropriate differential functions
 (see below for the latter notions). E. g., the
 Rellich's Identity for a function \cite{rel} or for a pair of
 functions \cite{em0,em} can be obtained in this way using a
 well-known Lagrangian and the radial vector field determining a
 dilation in $\mathbb{R}^n$. That is, both the Pohozaev's and the
 Rellich's identities come from the Noether's Identity. We claim that
 there is a certain kind of interplay between the integral identities of
 Pohozaev-Rellich type, Hardy-Sobolev
 Inequalities, Liouville type theorems, existence of
 conformal Killing vector fields and divergence symmetries of nonlinear Poisson
 equations on Riemannian manifolds. These connections will be
 studied in more details elsewhere.

 The present paper is a natural continuation of \cite{ye2}.
 Nevertheless, it can be read independently of \cite{ye2}.

 To begin with, let $M$ be an oriented Riemannian manifold of
 dimension $n\geq 3$ endowed with a metric $g=(g_{ij})$. We assume that
$M$ has a boundary $\kd M$ of class $C^{\infty}$. The local
coordinates of $M$ will be denoted by $x=(x^1,...,x^n)$. We denote
 by $dV$ and $dS$ the volume and surface measures
 with respect to the metric $g$, and by $\nu $ -
 the outward unit vector normal to $\kd M$.

 Now we introduce further notation and state the Noether's
 Identity which is the main ingredient of the proposed method.
 For more details see \cite{i1,i2,ol}.

 We consider a collection of smooth functions
 $u^{\al }(x)$, $\al =1,2,...,m$, defined on the manifold $M$. For an integer number $k\geq 1$ we
 let $u_{(k)}$ denote the set of all partial derivatives of $u^{\al }(x)$, up to order $k$.
 That is,
 \[ u_{(k)} = \left \{ u^{\al }_{i_1 ...i_s} \;|\; \al =1,2,...,m,
 \;s=1,...,k,\;i_1, ...,i_s = 1,...,n \right \}, \]
 where
 \[u^{\al }_{i_1 ...i_s}=\ds{\frac{{\kd }^s u^{\al }}{\kd x^{i_1}...\kd x^{i_s}}} \]
 and $u=u_{(0)}=(u^1,...,u^m)$.

 Following Olver \cite{ol} we introduce the notion of differential
 function.

\begin{defn} A smooth function of $x$, $u$ and derivatives
 of $u$ up to some finite, but unspecified order, is called
 differential function.
 \end{defn}

 The vector space of all differential functions of all orders is
 denoted by $\cal A$.

 We recall that the total derivative operator
 \[D_i =\frac{\kd}{\kd x_i}+u_i^{\al }\;\frac{\kd}{\kd
 u^{\al } } +u_{ij}^{\al } \;\frac{\kd}{\kd u_j^{\al }}+...+u_{i i_1 i_2 ... i_l}^{\al }\;
 \frac{\kd}{\kd u_{i_1 i_2 ... i_l}^{\al }} +...\;\; \]
 Above and throughout this paper we shall
suppose summation
 from $1$ to $n$ over repeated Latin indices and from $1$ to $m$
 over repeated Greek indices.

 Given $n$ differential functions ${\xi }^i = {\xi }^i (x,u, u_{(1)},...) \in \A$ and
 $m$ differential functions ${\eta }^{\al }={\eta }^{\al } (x,u, u_{(1)},...) \in \A$, let
 \[\begin{array}{lll}{\eta }^{(1)\al }_{i}&=&D_i{\eta }^{\al } -
 (D_i{\xi }^j)u^{\al }_j, \;\;i=1,2, ...,n;\\ & & \\
 {\eta }^{(l)\al }_{i_1 i_2...i_l} &=&
 D_{i_l} {\eta }^{(l-1)\al }_{i_1 i_2...i_{l-1}}
 - (D_{i_l}{\xi }^{j})u^{\al }_{i_1 i_2
 ...i_{l-1}j}\\ & & \\ &=& D_{i_1}D_{ i_2}...D_{ i_l} Q^{\al }
+ {\xi }^i u^{\al }_{i_1 i_2 ...i_k i},\end{array}\]
 where $i_l=1,2,...,n$ for $l=2,3,...,k$, $k=2,3, ...$ and
 $Q^{\al }={\eta }^{\al } - {\xi }^i u^{\al }_i $ are the Lie characteristic functions.

 Further one associates to ${\xi }^i$ and ${\eta }^{\al }$ the following
 partial differential operators acting on $\A$:

 - The operator:
 \[X={\xi }^i\frac{\kd }{\kd x_i}+ {\eta }^{\al }  \frac{\kd }{\kd u^{\al }}+
 \sum_{s=1}^{\infty }{\eta }^{(s)\al }_{i_1 i_2...i_s}
 \frac{\kd }{\kd u^{\al }_{i_1 i_2...i_s}}; \]

 - The Euler operator $E=(E_1,...,E_m) $, where $E_{\al }$ is defined by
 \[E_{\al } = \frac{\kd }{\kd u^{\al }} + \sum_{s=1}^{\infty } (-1)^s D_{i_1}D_{i_2}...D_{i_s}
 \frac{\kd }{\kd u^{\al }_{i_1 i_2...i_s}} ;\]

 - The Noether operator $N=(N^1,...,N^n) $, where:
 \[ N^i ={\xi }^i + Q^{\al }\left [ \frac{\kd }{\kd u^{\al }_i}
 + \sum_{s=1}^{\infty } (-1)^s D_{i_1}D_{i_2}...D_{i_s}
 \frac{\kd }{\kd u^{\al }_{i_1 i_2...i_s i}} \right ] \]
  \[+\sum_{r=1}^{\infty }D_{j_1}D_{j_2}...D_{j_r}Q^{\al }
 \left [ \frac{\kd }{\kd u^{\al }_{j_1 j_2...j_r i}}\right. \] \[\left.  + \sum_{s=1}^{\infty }
 (-1)^s D_{l_1}D_{l_2}...D_{l_s}
 \frac{\kd }{\kd u^{\al }_{l_1 l_2...l_s j_1 j_2... j_r i}} \right ]. \] 

 These operators are related by the Noether's Identity
 \cite{i1,i2}:
 \bb\label{n1} X+D_i{\xi }^i=Q^{\al }E_{\al }+D_i N^i .\ee

 The identity (\ref{n1}) was explicitly stated for the first time in the work \cite{i1} by
 Ibragimov who named it in honour of E. Noether. As Ibragimov has pointed out it is clear that
 this identity makes the proof of the Noether Theorem \cite{n,ol} purely
 algebraic and very simple.

 Now we outline the proposed Noetherian approach to Rellich's
 Identity. Let $u$ be a vector-valued smooth function and $P=(P_1,...,P_m)$ - a
 linear or nonlinear partial differential operator of an arbitrary
 even order $2k$. Let $L\in \A$ be a differential function of
 order $k$ such that
 \[Pu = E(L), \]
 where $E$ is the Euler operator (see above). The function $L$
 will play the role of a `Lagrangian'. Then one writes the
 Noether's Identity
 corresponding to $L$ and a suitable operator $X$
 \[ XL+L\;D_i{\xi }^i = ({\eta }^{\al } - {\xi }^i u^{\al }_i) P_{\al }u +D_i(N^i L),\]
 integrates and applies the divergence theorem:
 \bb\label{a22} 
 \ds{\int_M} [ XL+(D_i{\xi }^i)L ]dV = \ds{\int_M}({\eta }^{\al } - {\xi }^i u^{\al }_i) P_{\al }u\;dV +\ds{\int_{\kd M}} (N^i L){\nu }_i dS. 
 \ee 
 
 The identity (\ref{a22}) is the most general form of the Rellich's Identity. We emphasize that in (\ref{a22}) one has the freedom to chose the operator $X$ (that is, the differential functions ${\xi }^i$ and ${\eta }^{\al }$) as well as the Lagrangian $L$ depending on the specific research objective. In this way various Rellich type identities can be obtained, one of which is established in the following  
 
\begin{thm} Let $u,v\in C^2(\bar{M})$ be two given
 functions, $h =h^i(x) \frac{\kd }{\kd x^i}$ a $C^1(\bar{M})$ vector field, 
 and $F=F(x, u_{(1)}, v_{(1)})\in \cal{A}$. Then the
 following  identity holds:  
 \begin{eqnarray} \label{n9}
 \ds{\int_M} \{ div F_p\;(h,F_q) dV &+& div F_q\;(h,F_p)\} dV =
 -\ds{\int_M} {\cal L}_h g_{ik}\; F_{u_i} F_{v_j} dV\nonumber \\
 & & \nonumber \\
 & - & \ds{\int_M} h_j \{ F_{u_i} {\nb }_i F_{v_j} + F_{v_i} {\nb }_i F_{u_j} \}
 dV  \\
 & & \nonumber \\
 &+& \ds{\int_{\kd M}} \{ (F_p,\nu ) (h,F_q)+(F_p,\nu ) (h,F_q)
 \}dS,
 \nonumber\end{eqnarray}
 {\it where } \[ F_p = \left ( \frac{\kd F}{\kd u_1}(x,u_1,...u_n),..., \frac{\kd F}{\kd u_n}(x,u_1,...u_n) \right ), \] \[ F_q = \left ( \frac{\kd F}{\kd v_1}(x,v_1,...v_n),..., \frac{\kd F}{\kd v_n}(x,u_1,...u_n) \right ),\] ${\nabla }^i$ {\it is the covariant derivative corresponding to
the Levi-Civita connection, uniquely determined by $g$, and ${\cal L}_h g_{ij} $ is the Lie derivative of the metric $g_{ij}$ 
 with respect to the vector field $h$.}
 \end{thm}
 
In the Euclidean case, the identity (\ref{n9}) coincides with the main identity in \cite{em}, up to an integration by parts. 
 
One more time we observe that the smooth vector field $h$ can be chosen in an arbitrary way which provides a wide variety of possibilities for applications, even in the Euclidean case!
 
However, we shall not proceed further in this generality. Our next purpose is to make use of conformal Killing vector fields in order to obtain Riemannian analogs of some Rellich type integral identities established in \cite{em}. \vspace{-0.4cm}
 
We recall that a nonisometric conformal Killing
vector field on $M$ is a vector field $h=h^i\frac{\kd}{\kd x^i}$ which satisfies 
\bb\label{ww} 
{\nb }^i h^j + {\nb }^j h^i= \frac{2}{n} (div (h))\;g^{ij} =\mu\; g^{ij},\; \mu\neq 0.  
\ee
Here $div(h)= {\nb }_i h^i$ is the covariant divergence operator and $g^{ij}$ is the inverse matrix of $g_{ij}$.

We would like to point out that, in fact, the conformal Killing vector field generalizes the dilations in ${\mathbb{R}}^n$, determined by $x^i\frac{\kd }{\kd x^i}$, and enables us to follow in the general manifold case some of the arguments successfully used in the Euclidean case. To show this, for example, we obtain a higher order Rellich type identity involving the
polyharmonic operator on Riemannian manifolds admitting proper homothetic
 transformations, that is, transformations whose infinitesimal
 generators $h$ satisfy (\ref{ww}) with $\mu = const \neq 0$. It
 is easy to see that in this case one can assume $div(h)=n$
 without loss of generality. (Note that $h=x^i\frac{\kd }{\kd x^i}=r\frac{\kd }{\kd r}$,  corresponding to a dilational transformation in ${\mathbb{R}}^n$, satisfies the last relation.) 
 
\begin{thm} Let $u,v\in C^{4m}(\bar{M})$ be two given
 functions and let  $h =h^i(x) \frac{\kd }{\kd x^i}$ be a $C^1(\bar{M})$
 contravariant conformal Killing vector field such that $div(h)=n$.

Then the following  identity holds:
 \begin{eqnarray}\label{n30}
 R_{2m}(u,v)&=&(4m-n)\ds{\int_M} {\De }_g^{2m}
 u\;v\;dV+ \ds{\int_{\kd M}} \lm u\lm v (h,\nu )\;dS\nonumber \\
 & & \nonumber \\
 & + &(4m-n)\sum_{l=0}^{m-1}\ds{\int_{\kd M}}{\De }_g^{m+l}u\; (\nb
 ({\De }_g^{m-l-1} v),\nu )\;dS\nonumber \\
 & & \nonumber \\
 & - & (4m-n)\sum_{l=0}^{m-1}\ds{\int_{\kd M}}{\De }_g^{m-l-1}v\; (\nb
 ({\De }_g^{m+l} u),\nu )\;dS \nonumber \\
 & & \nonumber \\
 &+&\sum_{l=0}^{m-1}\ds{\int_{\kd M}}(2l\; {\De }_g^l v + h^k{\nb }_k({\De
 }_g^l v))\;({\nb }({\De }^{2m-1-l}_g
 u),\nu )\;dS  \\
 & & \nonumber \\
 &-&\sum_{l=0}^{m-1}\ds{\int_{\kd M}}(2l\;({\nb } ({\De }^{l}_g v),\nu ) +{\nb }^i h^k {\nb
 }_k ({\De }^{l}_g v){\nu }_i \nonumber \\ & &
 \;\;\;\;\;\;\;\;\;\;\;\;\;\;\;\;\;\;\;\;\;\;\;\;\;\;\;\;\;\;\;\;\;\;\;\;\;+h^k {\nb }^i{\nb }_k ({\De }^{l}_g
 v){\nu }_i)
 \;({\De }^{2m-1-l}_g u)\;dS\nonumber \\
 &+&\sum_{l=0}^{m-1}\ds{\int_{\kd M}}(2l\; {\De }_g^l u + h^k{\nb }_k({\De
 }_g^l u))\;({\nb }({\De }^{2m-1-l}_g
 v),\nu )\;dS \nonumber \\
 & & \nonumber \\
 &-&\sum_{l=0}^{m-1}\ds{\int_{\kd M}}(2l\;({\nb } ({\De }^{l}_g u),\nu ) +{\nb }^i h^k {\nb
 }_k ({\De }^{l}_g u){\nu }_i \nonumber \\ & &
 \;\;\;\;\;\;\;\;\;\;\;\;\;\;\;\;\;\;\;\;\;\;\;\;\;\;\;\;\;\;\;\;\;\;\;\;\;+ h^k {\nb }^i{\nb }_k ({\De }^{l}_g
 u){\nu }_i)
 \;({\De }^{2m-1-l}_g v)\;dS.\nonumber\end{eqnarray}
 {\it where }
 \[R_{2m}(u,v)=\ds{\int_M} \{ {\De }_g^{2m} u\;(h,\nb v) +
 {\De }_g^{2m} v\;(h,\nb u)\} dV \]
 {\it and ${\De }_g$ is the Laplace-Beltrami operator corresponding to the metric $g$.}
 \end{thm}
 
A similar identity involving odd powers of ${\De }_g$ is also
 valid. We omit it in order not to increase the volume of this
 paper.

Let $u=v$. From (\ref{n30}) we immediately obtain
 
\begin{cor} Let $u\in C^{4m}(\bar{M})$ be a given
 function and let $h =h^i(x) \frac{\kd }{\kd x^i}$ be a $C^1(\bar{M})$
 contravariant conformal Killing vector field such that $div(h)=n$.

 Then the following  identity holds:
 \[ \ds{\int_M} {\De }_g^{2m}u(h,\nb u)\;dV =\frac{4m-n}{2}\ds{\int_M} {\De }_g^{2m}
 u\;u\;dV + \frac{1}{2}\ds{\int_{\kd M}} (\lm u )^2 (h,\nu )\;dS\]
 \begin{eqnarray}\label{n31}
 &+& \frac{4m-n}{2}\sum_{l=0}^{m-1}\ds{\int_{\kd M}}{\De }_g^{m+l}u\; (\nb
 ({\De }_g^{m-l-1} u),\nu )\;dS
 \nonumber \\
 & & \nonumber \\
 & - & \frac{4m-n}{2}\sum_{l=0}^{m-1}\ds{\int_{\kd M}}{\De }_g^{m-l-1}u\; (\nb
 ({\De }_g^{m+l} u),\nu )\;dS \nonumber \\
 & &  \\
 &+& \sum_{l=0}^{m-1}\ds{\int_{\kd M}}(2l\; {\De }_g^l u + h^k{\nb }_k({\De
 }_g^l u))\;({\nb }({\De }^{2m-1-l}_g
 u),\nu )\;dS \nonumber \\
 & & \nonumber\\
 &-&\sum_{l=0}^{m-1}\ds{\int_{\kd M}}(2l\;({\nb } ({\De }^{l}_g u),\nu ) +{\nb }^i h^k {\nb
 }_k ({\De }^{l}_g u){\nu }_i\nonumber \\
& & \nonumber \\
 & &+ h^k {\nb }^i{\nb }_k ({\De }^{l}_g
 u){\nu }_i)
 \;({\De }^{2m-1-l}_g u)\;dS.\nonumber\end{eqnarray}
 \end{cor}
 
The identity (\ref{n31}) in the Euclidean case with $h=x^i\frac{\kd}{\kd x^i}$ is similar to an
 identity obtained in \cite{em} and used in \cite{em3} to
 establish a nonexistence result for positive radial solutions of
 semilinear polyharmonic equations in $\mathbb{R}^n$.

Further we prove a biharmonic Rellich identity in a more general context.

\begin{thm} Let $u,v\in C^{4}(\bar{M})$ be two given
 functions and let $h =h^i(x) \frac{\kd }{\kd x^i}$ be a
 $C^1(\bar{M})$ conformal Killing vector field satisfying (\ref{ww}). Then the
 following  identity holds:
 \[ \ds{\int_M} \{ {\De }_g^2 u\;(h,\nb v) dV + {\De }_g^2 v\;(h,\nb u)\} dV =
 \frac{4-n}{2}\ds{\int_M} \mu\;{\De }_g u{\De }_g v dV \]
 \begin{eqnarray}\label{n38}
 &+&\frac{1}{n-1}\ds{\int_M} ({\cal L}_h R+\mu R) (u{\De }_g v +v{\De }_g u)dV \nonumber
 \\
 & & \nonumber \\
 & -& \ds{\int_M}(u(\nb\mu , \nb ({\De }_g v)) + v(\nb\mu , \nb ({\De }_g
 u)))dV \nonumber \\
 & & \\
 &+& \ds{\int_{\kd M}}({\De }_g u{\De }_g v (h,\nu ) +
 (u{\De }_g v +v{\De }_g u)(\nb\mu ,\nu) \nonumber\\
 & & \nonumber\\
 &+& (h,\nb v) {\nb }^i({\De }_g u) {\nu }_i -{\De }_g u ({\nb }^i h^k. v_k{\nu }_i +h^k
 {\nb }_k{\nb }^i v .{\nu }_i)\nonumber\\
 & & \nonumber\\
 &+& (h,\nb u) {\nb }^i({\De }_g v) {\nu }_i -{\De }_g v ({\nb }^i h^k. u_k{\nu }_i +h^k
 {\nb }_k{\nb }^i u .{\nu }_i))dS, \nonumber \end{eqnarray}
 {\it where $R$ is the scalar curvature of $M$ and ${\cal L}_h R$ is its Lie derivative with respect to the vector field $h$.}
\end{thm}

Setting $u=v$ in (\ref{n38}) we immediately obtain

\begin{cor} {\it   Let $u\in C^4(\bar{M})$ be a given
 function and let $h =h^i(x) \frac{\kd }{\kd x^i}$ be a $C^1(\bar{M})$
 conformal Killing vector field satisfying (\ref{ww}). Then the
 following  identity holds:}
\[ \ds{\int_M} {\De }_g^2 u\;(h,\nb u) dV =\frac{4-n}{4}\ds{\int_M} \mu\;({\De }_g
  u)^2dV -\ds{\int_M} u(\nb\mu , \nb ({\De }_g u))dV \]
\begin{eqnarray}\label{n39}
 &+&\frac{1}{n-1}\ds{\int_M} ({\cal L}_h R+\mu R) u{\De }_g u\;dV \nonumber \\
 & & \\
 &+&\ds{\int_{\kd M}}(\frac{1}{2}({\De }_g u)^2 (h,\nu ) + u{\De }_g
 u (\nb\mu ,\nu)\nonumber\\
 & & \nonumber \\&+&(h,\nb u) {\nb }^i({\De }_g u) {\nu }_i -
 {\De }_g u ({\nb }^i h^k. u_k{\nu }_i +h^k {\nb }_k{\nb }^i u .{\nu }_i)) dS. \nonumber
 \end{eqnarray}
 \end{cor}

This paper is organized as follows. In section 2 we prove Theorem 1.2. Then, aiming to show how the proposed approach works,
 we re-establish in section 3 two known integral identities
obtained in \cite{ye2} and \cite{em}. The proofs of theorems 1.3 and 1.5 are presented in sections 4 and 5 respectively. In section 6 we establish a nonexistence result for systems of two biharmonic 
equations on Riemannian manifolds.
 
For further applications, e. g. Hardy and Caffarelli-Kohn-Nirenberg type inequalities, see \cite{ye2} and \cite{ydb1}.
 
\section{A general Rellich type identity}
 
In this section we prove Theorem 1.2. 
 
In the Noether's identity (\ref{n1}) we take $m=2$, $u^1=u$ and
 $u^2=v$. Consider the vector field
 \bb\label{n2} X ={\xi }^i(x) \frac{\kd }{\kd x_i}+ \eta
 (x,u,v,u_{(1)}, v_{(1)}) \frac{\kd }{\kd u} + \phi (x,u,v,u_{(1)}, v_{(1)}) \frac{\kd }{\kd
 v} \ee
 where ${\xi }^i,{\eta }^1=\eta ,{\eta }^2=\phi\in \A  $. Then the Noether's Identity
 (\ref{n1}) applied to $ L=L(x,u,v,u_{(1)}, v_{(1)})$, another
 differential function, reads:
 \bb\label{n3}\begin{array}{lll} X^{(1)}L +L D_i{\xi }^i &= & E_u (L) (\eta - {\xi }^j u_j) +
 E_v (L) (\phi - {\xi }^j u_j)\\ & & \\ & +&D_ i \ds{\left [{\xi }^i L +
 \frac{\kd L}{\kd u_i} (\eta - {\xi }^j u_j) +
 \frac{\kd L}{\kd v_i} (\phi - {\xi }^j v_j) \right ]}, \end{array}\ee
 where
 \[ D_i =\frac{\kd}{\kd x_i}+u_i\;\frac{\kd}{\kd u } +v_i\;\frac{\kd}{\kd v }
 +u_{ij}\;\frac{\kd}{\kd u_j}+v_{ij}\;\frac{\kd}{\kd v_j}...\] 
 \[+
 u_{i i_1 i_2 ... i_l}\; \frac{\kd}{\kd u_{i_1 i_2 ... i_l}} +
 v_{i i_1 i_2 ... i_l}\; \frac{\kd}{\kd v_{i_1 i_2 ... i_l}} +... \]
 is the total derivative operator, $E=(E_u,E_v)$ is the Euler
 operator with components
 \[E_{u } = \frac{\kd }{\kd u} + \sum_{s=1}^{\infty } (-1)^s D_{i_1}D_{i_2}...D_{i_s}
 \frac{\kd }{\kd u_{i_1 i_2...i_s}} \]
 and
 \[E_{v} = \frac{\kd }{\kd v} + \sum_{s=1}^{\infty } (-1)^s D_{i_1}D_{i_2}...D_{i_s}
 \frac{\kd }{\kd v_{i_1 i_2...i_s}}, \]
 and $X^{(1)}$ is the first order prolongation of $X$ given by
 \[ X^{(1)}=X + {\eta }^{(1)}_{i} \frac{\kd }{\kd u_i}
 + {\phi }^{(1)}_{i} \frac{\kd }{\kd v_i}.\]

 Now we shall choose special forms of the differential functions
 involved in (\ref{n3}). Namely, let $F=F(x,u_{(1)}, v_{(1)})\in\A $ and
 $h =h^i(x) \frac{\kd }{\kd x^i} $ be a smooth vector filed.
 Then we set $L=F$, ${\xi }^i=0 $,
 \[ \eta = (h,F_q) = h_j F_{v_j}=g_{ij} h^i F_{v_j}, \;\;\;
 \phi = (h,F_p) = h_j F_{u_j}=g_{ij} h^i F_{u_j}.\]
 The differential operator $X$ assumes the form,
 \[ X=(h,F_q)\frac{\kd }{\kd u} + (h,F_p)\frac{\kd }{\kd v}, \]
 and the General Prolongation Formula (\cite{ol}, p. 113) immediately
 gives its first order prolongation:
 \[X^{(1)}=X + D_i (h,F_q)\;\frac{\kd }{\kd u_i}
 + D_i (h,F_p)\; \frac{\kd }{\kd v_i}. \]
 We substitute this data into the Noether's Identity (\ref{n3}).
 The left-hand side of (\ref{n3}) is
 \begin{eqnarray}\label{n4}X^{(1)}L +L D_i{\xi }^i&=&D_i (h,F_q) F_{u_i} + D_i (h,F_p) F_{v_i}\nonumber \\ 
 & & \\ & = & (F_p, \nb (h,F_q)) + (F_q, \nb (h,F_p))\nonumber \end{eqnarray}
 while its right-hand side is given by,
 \bb\label{n5} \begin{array}{lll}
 E_u (L)(h,F_q)&+ &E_v (L)(h,F_p)+D_i [F_{u_i} (h,F_q)+ F_{v_i} (h,F_p)]\\ & & \\
 & =& - D_i F_{u_i}(h,F_q)-  D_i F_{v_i}(h,F_p)\\ & & \\ & &+
 D_i [F_{u_i} (h,F_q)+ F_{v_i} (h,F_p) ] \\ & & \\ & =&
 -div\; F_p\; (h,F_q) -div\; F_q \;(h,F_p)\\ & & \\ & &+D_i [F_{u_i} (h,F_q)+ F_{v_i} (h,F_p)].\end{array}\ee
 From (\ref{n3}), (\ref{n4}) and (\ref{n5}) it follows that
 \bb\label{n6}\begin{array}{lll} div\;F_p \;(h,F_q) +div\;F_q\; (h,F_p)&= &
 - (F_p, \nb (h,F_q)) - (F_q, \nb (h,F_p))\\ & & \\
 &+&D_i [F_{u_i} (h,F_q)+ F_{v_i} (h,F_p)].\end{array}\ee
 Actually, the identity (\ref{n6}) itself is obvious! Nevertheless
 its usefulness comes from the fact that it is directly obtained
 from the Noether's Identity, providing in this way a general
 procedure to get Rellich type identities, as it can be seen
 below.

 We observe that the identity (\ref{n6}) holds if we replace the
 partial derivatives $D_i$ by the covariant derivatives ${\nb }_i$
 corresponding to the Levi-Civita connection of $M$. In fact, we
 could have done this in the Noether's Identity from the beginning
 of this procedure noting that the function $L$ contains only the
 first derivatives of $u$ and $v$. We shall come back to this point
 later.

 Differentiating in (\ref{n6}) and interchanging some of the
 indices we obtain: 
 \begin{eqnarray} \label{n8}
  div F_p\;(h,F_q) dV &+& div F_q\;(h,F_p) =
 - \{ {\nb }_i h_j + {\nb }_j h_i \}F_{u_i} F_{v_j} \nonumber \\
 & & \nonumber \\
 & - &  h_j \{ F_{u_i} {\nb }_i F_{v_j} + F_{v_i} {\nb }_i F_{u_j} \}
   \\
 & & \nonumber \\
 &+&  {D }_i [F_{u_i} (h,F_q)+ F_{v_i} (h,F_p) ] .
 \nonumber\end{eqnarray}
 We observe that in (\ref{n8}) the Lie derivative of the metric
 with respect to the vector field $h$ appears, namely
 \[ {\cal L}_h g_{ij} = {\nb }_i h_j + {\nb }_j h_i.
 \]
 Thus
 \begin{eqnarray} \label{ab8}
  div F_p\;(h,F_q) dV &+& div F_q\;(h,F_p) =
 - {\cal L}_h g_{ij}\;F_{u_i} F_{v_j} \nonumber \\
 & & \nonumber \\
 & - &  h_j \{ F_{u_i} {\nb }_i F_{v_j} + F_{v_i} {\nb }_i F_{u_j} \}
   \\
 & & \nonumber \\
 &+&  {D }_i [F_{u_i} (h,F_q)+ F_{v_i} (h,F_p) ] .
 \nonumber\end{eqnarray}
 
Then by integrating the identity (\ref{ab8}) and applying the divergence theorem, we get 
(\ref{n9}).

\section{A Rellich type identity involving the Laplace operator}
 
In this section we illustrate the proposed Noetherian approach to
 Rellich type identities recovering two known integral identities
 obtained in \cite{ye2} and \cite{em}.

 For this purpose let us now suppose that $h$ in the preceding section is a conformal Killing vector field:
 \[ {\cal L}_h g_{ij} ={\nb }_i h_j + {\nb }_j h_i= \frac{2}{n}
 (div\; h)g_{ij}. \]
 Hence and from (\ref{n9}):
 \begin{eqnarray} \label{n10}
 \ds{\int_M} \{ div F_p\;(h,F_q) dV &+& div F_q\;(h,F_p)\} dV =
 -\ds{\frac{2}{n}}\ds{\int_M} div\;h\;g_{ij} F_{u_i} F_{v_j} dV\nonumber \\
 & & \nonumber \\
 & - & \ds{\int_M} h_j \{ F_{u_i} {\nb }_i F_{v_j} + F_{v_i} {\nb }_i F_{u_j} \}
 dV  \\
 & & \nonumber \\
 &+& \ds{\int_{\kd M}} \{ (F_p,\nu ) (h,F_q)+(F_p,\nu ) (h,F_q)
 \}dS.
 \nonumber\end{eqnarray}
 Let  \[ F= g^{ij}(u_i u_j +v_i v_j)/2. \]
 Then
 \[ F_{u_i}=g^{is} u_s=u^i,\;\;\; F_{v_k}=g^{ks} v_s=v^k,\;\;\;
 (h,F_p)=h^k u_k,\;\;\; (h,F_q)=h^k v_k, \]
 \[ div F_p = {\nb }_i u^i={\De }_g u, \;\;\;\;div F_q = {\nb }_i v^i={\De }_g
 v. \]
 With this choice of $F$ the identity (\ref{n10}) yields
 \begin{eqnarray} \label{n11}
 \ds{\int_M} \{ {\De }_g u\;(h,\nb v) dV &+& {\De }_g v\;(h,\nb u)\} dV =
 -\ds{\frac{2}{n}}\ds{\int_M} div\;h\;(\nb u,\nb v) dV\nonumber \\
 & & \nonumber \\
 & - & \ds{\int_M} h_j (u^i {\nb }_i v^j + v^i {\nb }_i u^j)
 dV  \\
 & & \nonumber \\
 &+& \ds{\int_{\kd M}} \{ (F_p,\nu ) (h,F_q)+(F_p,\nu ) (h,F_q)
 \}dS.
 \nonumber\end{eqnarray}
 Further we integrate by parts in the term of the second line of
 (\ref{n11}) taking into account the fact that the second covariant
 derivatives of a function commute:
 \begin{eqnarray} \label{n12} - \ds{\int_M} h_j (u^i {\nb }_i v^j &+& v^i {\nb }_i u^j)
 dV = -\ds{\int_M} h^k {\nb }_k (\nb u,\nb v) dV \nonumber \\
 & & \\
 &=&\ds{\int_M} div\;h\;(\nb u,\nb v) dV -\ds{\int_{\kd M}} (\nb u,\nb
 v)(h,\nu ) dS. \nonumber\end{eqnarray}
 Substituting (\ref{n12}) into (\ref{n11}) we obtain
 \[\!\!\!\!\!\!\!\!\!\!\!\!\!\!\!\!\!\!\ds{\int_M} \{ {\De }_g u\;(h,\nb v) dV + {\De }_g v\;(h,\nb u)\} dV \]
 \begin{eqnarray}\label{n13}
 & & \nonumber \\
 &=& \frac{n-2}{n}\ds{\int_M} div\;h\;(\nb u,\nb v) dV 
 \\
 & & \nonumber \\
 &+& \ds{\int_{\kd M}} \{ \frac{\kd u}{\kd\nu} (h,\nb v)+
 \frac{\kd v}{\kd\nu}(h,\nb u) -(\nb u,\nb v)(h,\nu )
 \}dS. \nonumber \end{eqnarray}

 Clearly (\ref{n13}) is the identity (17) of \cite{ye2}. Moreover,
 if $M=\Om $ is a bounded
 domain in ${\mathbb{R}}^n$, $g_{ij}={\de }_{ij}-$the Euclidean
 metric and $h =x^i \frac{\kd }{\kd x^i} $, then we obtain the Rellich type identity
 (2.5) established in \cite{em}, pp. 128-129.

 Before concluding this section we would like to observe that in
 \cite{ye2} the identity (\ref{n13}) is obtained following the
 argument in \cite{em} while here it is obtained by applying the
 proposed Noetherian approach to Rellich identities.

\section{A higher order Rellich type identity}

In this section we prove Theorem 1.3, namely, we obtain a Rellich type identity involving the
 polyharmonic operator ${\De }^k_g $ on a Riemannian manifold
 $(M^n,g)$, where $k\geq 2$ is an even number. The case $k-$odd
 can be treated in a similar way.

 As in the preceding section, we take $u^1=u$, $u^2=v$, ${\xi
 }^i=0$, ${\eta }^1=\eta = (h,\nb v)$, ${\eta }^2=\phi = (h,\nb v)$
 in the Noether's Identity (\ref{n1}) which is applied to the
 differential function
 \[ L= \frac{1}{2}(\lm u)^2 + \frac{1}{2}(\lm v)^2,\]
 where $m=k/2\geq 1$. Then the Noether's Identity reads
 \bb\label{n14} X^{(k)}L=E_u(L)\eta +E_v(L)\phi + {\nb }_i w^i, \ee
 where $ X^{(k)}$ is the $k-$th order prolongation of \[ X=\eta
 \frac{\kd }{\kd u}+\phi\frac{\kd }{\kd v}, \]
 $(E_u,E_v) $ is the Euler operator and $w^i=N^i (L)$, ($N=(N^1,...,N^n)$ being
 the Noether operator) will be explicitly calculated later.

 Calculations similar to those presented in \cite{ydb} lead to
 \bb\label{n15}X^{(k)}L=(\lm u)(\lm \eta )+ (\lm v)(\lm \phi ) \ee
 and
 \bb\label{n16}E_u(L)\eta +E_v(L)\phi = {\De }^{2m}_g u.\eta +{\De }^{2m}_g
 v.\phi . \ee
 Thus (\ref{n14}), (\ref{n15}) and (\ref{n16}) imply the following.

\begin{prop}{\it If $u,v\in C^{4m}(M)$, then }
 \bb\label{n17} (\lm u)(\lm \eta )+ (\lm v)(\lm \phi )=
 {\De }^{2m}_g u.\eta +{\De }^{2m}_g
 v.\phi + {\nb }_i w^i, \ee
 {\it where} $\eta = (h,\nb v)$, $\phi = (h,\nb v)$ {\it and}
 \begin{eqnarray}\label{n18}
 w^i & = & -\sum_{s=0}^{m-1}{\De }^s_g\eta \;{\nb }^i ({\De }^{2m-1-s}_g
 u) +\sum_{s=0}^{m-1}{\nb }^i({\De }^s_g\eta )\;{\De }^{2m-1-s}_g u \nonumber \\
 & & \\
 &-&\sum_{s=0}^{m-1}{\De }^s_g\phi \;{\nb }^i ({\De }^{2m-1-s}_g
 v) +\sum_{s=0}^{m-1}{\nb }^i({\De }^s_g\phi )\;{\De }^{2m-1-s}_g v. \nonumber
 \end{eqnarray} 
 \end{prop}

Now we would like to comment on an important point. To write
 (\ref{n17}) we have substituted in the Noether's Identity
 (\ref{n1}) the partial derivatives $D_i$ by the covariant
 derivatives ${\nb }_i$. In this way we have actually used a {\it
 covariant} Noether's Identity which up to our knowledge has not
 been previously established and used in the literature. Therefore it
 requires a rigorous proof. However such a proof is very lengthy
 since it contains a lot of technical details, in particular,
 details regarding its application to the polyharmonic Lagrangian
 $L$ (as in \cite{ydb}) as well as careful commutations of the
 covariant derivatives which appear during this procedure. For
 these reasons we shall not present such a proof here merely
 pointing out that such a problem (viz. study of covariant Noether's Identity)
 will be treated elsewhere. Nevertheless, with (\ref{n17}) at
 hand, we can prove it directly using a simple alternative
 argument.
 
{\bf Proof of Proposition 4.1.} We calculate the divergence
 \begin{eqnarray}
 {\nb }_ i w^i  = & -&\sum_{s=0}^{m-1}{\nb }_ i{\De }^s_g\eta \;{\nb }^i ({\De }^{2m-1-s}_g
 u) -\sum_{s=0}^{m-1}{\De }^s_g\eta \; ({\De }^{2m-s}_g
 u)\nonumber \\
 & & \nonumber \\
 &+& \sum_{s=0}^{m-1}{\De }^{s+1}_g\eta \;{\De }^{2m-1-s}_g u
 +\sum_{s=0}^{m-1}{\nb }^i{\De }^{s}_g\eta \;{\nb }_i({\De }^{2m-1-s}_g u) \nonumber \\
 & & \nonumber \\
 &-&\sum_{s=0}^{m-1}{\nb }_ i{\De }^s_g\phi \;{\nb }^i ({\De }^{2m-1-s}_g
 v) -\sum_{s=0}^{m-1}{\De }^s_g\phi \; ({\De }^{2m-s}_g
 v)\nonumber \\
 & & \nonumber \\
 &+& \sum_{s=0}^{m-1}{\De }^{s+1}_g\phi \;{\De }^{2m-1-s}_g v
 +\sum_{s=0}^{m-1}{\nb }^i{\De }^{s}_g\phi \;{\nb }_i({\De }^{2m-1-s}_g v) \nonumber \\
 =& -&\sum_{s=1}^{m-1}{\De }^s_g\eta \; ({\De }^{2m-s}_g
 u) - \eta {\De }^{2m}_g u + \sum_{s=1}^{m}{\De }^{s}_g\eta \;{\De }^{2m-s}_g u\nonumber \\
 & & \nonumber \\
 & -&\sum_{s=1}^{m-1}{\De }^s_g\phi \; ({\De }^{2m-s}_g
 v) - \phi {\De }^{2m}_g v + \sum_{s=1}^{m}{\De }^{s}_g\phi \;{\De }^{2m-s}_g
 v\nonumber \\
 & & \nonumber \\
 =& -&\sum_{s=1}^{m-1}{\De }^s_g\eta \; ({\De }^{2m-s}_g
 u) - \eta {\De }^{2m}_g u + \sum_{s=1}^{m-1}{\De }^{s}_g\eta \;{\De }^{2m-s}_g
 u + \lm\eta\lm u\nonumber \\
 & & \nonumber \\
 & -&\sum_{s=1}^{m-1}{\De }^s_g\phi \; ({\De }^{2m-s}_g
 v) - \phi {\De }^{2m}_g v + \sum_{s=1}^{m-1}{\De }^{s}_g\phi \;{\De }^{2m-s}_g
 v + \lm\phi\lm v\nonumber \\
 & & \nonumber \\
 =&-& \eta {\De }^{2m}_g u - \phi {\De }^{2m}_g v + \lm\eta\lm u + \lm\phi\lm v,\nonumber \end{eqnarray}
 which completes the proof.

The next step is to calculate $\lm \eta $ and $\lm \phi $ appearing in (\ref{n17}). This is done in a sequence of lemmas and propositions.

\begin{lem} {\it If $h$ is a conformal Killing vector field
 satisfying} \bb\label{n19} \nabla^{i}h^{k}+\nabla^{k}h^{i}=\mu
g^{ij}=\frac{2}{n}\;div(h)\; g^{ik}, \ee
  {\it then }\bb\label{d0}{\De }_g
 h^{k}+R^{k}_{\;s}h^{s}=\frac{2-n}{2}g^{kj}\mu_{j}, \ee
 {\it where $R^{i}_{\;j} $ is the Ricci tensor}.
 \end{lem}

 \begin{proof} See \cite{yi}.
 \end{proof}

In particular, if the function $\mu $ = constant, we have
 \bb\label{n20}{\De }_g
 h^{k}=-R^{k}_{\;s}h^{s}.\ee

\begin{lem} {\it If $h$ is a conformal Killing vector field
 satisfying (\ref{n19}), then for any $v\in C^{2l+2}(M)$} we have,
 \bb\label{n21} 2{\nb }^i h^k\; {\nb }_i {\nb }_k ({\De }_g^l v)=\mu\; {\De }_g^{l+1} v.
 \ee
\end{lem}

 \begin{proof} The equality (\ref{n21}) is obtained by
  multiplying (\ref{n19}) by ${\nb }_i {\nb }_k ({\De }_g^l v) $
 and changing some of the indices.
\end{proof}

\begin{lem} {\it For any $v\in C^{2l+3}(M)$ the following identity holds,}
 \bb\label{n22} {\De }_g ({\nb }_k ({\De }_g^l v)) = {\nb }_k ({\De }_g^{l+1} v)
 + R^{s}_{\;k} {\nb }_s ({\De }_g^l v).\ee
\end{lem}

\begin{proof} The equality (\ref{n22}) for $l=0$ follows from the well-known
 commutation relation
 \[ (\nabla_{k}\nabla_{l}-\nabla_{l}\nabla_{k})T_{i}=-R_{i\;kl}^{\;s}
 T_{s}\]
 valid for any covariant field $T=(T_s)$. Here $ R_{i\;kl}^{\;s}$
 is the Riemann tensor. More generally, for a $p-$contravariant
 and $q-$covariant tensor
 \[ T=(T^{i_1...i_p}_{\;\;\;\;\;\;\;\;\;\;j_1...j_q})\]
 we have the commutation formula
 \begin{eqnarray}\label{d1}
 {\nb }_i{\nb }_k T^{i_1...i_p}_{\;\;\;\;\;\;\;\;\;\;j_1...j_q} & = &
 {\nb }_k{\nb }_i T^{i_1...i_p}_{\;\;\;\;\;\;\;\;\;\;j_1...j_q}+
 T^{s i_2...i_p}_{\;\;\;\;\;\;\;\;\;\;j_1...j_q}R^{i_1}_{\;\;ski}+...\nonumber \\
 & & \nonumber \\
 &+&T^{i_1...i_{p-1}s}_{\;\;\;\;\;\;\;\;\;\;\;\;\;\;\;j_1...j_q}R^{i_p}_{\;\;ski}
 - T^{i_1...i_p}_{\;\;\;\;\;\;\;\;\;\;s j_2...j_q}R^{s}_{\;\;j_1 k
 i} \\
 & & \nonumber \\
 &-&...-T^{i_1...i_p}_{\;\;\;\;\;\;\;\;\;\;j_1...j_{q-1}s}R^{s}_{\;\;j_q k
 i}.\nonumber \end{eqnarray}
 Setting in (\ref{d1}) $p=q=l$ and
 \[T^{i_1...i_l}_{\;\;\;\;\;\;\;\;\;\;j_1...j_l}={\nb }^{i_1}{\nb }_{j_1}
 ...{\nb }^{i_l}{\nb }_{j_l} v\]
 we obtain
  \begin{eqnarray}\label{d2}
 {\nb }_i{\nb }_k ({\nb }^{i_1}{\nb }_{j_1}
 ...{\nb }^{i_l}{\nb }_{j_l} v ) & =&{\nb }_k{\nb }_i({\nb }^{i_1}{\nb }_{j_1}
 ...{\nb }^{i_l}{\nb }_{j_l} v )\nonumber\\ & & \nonumber\\
 & +& {\nb }^s {\nb }_{j_1}({\nb }^{i_2}{\nb }_{j_2}
 ...{\nb }^{i_l}{\nb }_{j_l} v )R^{i_1}_{\;\;ski}+...\nonumber \\
 & & \nonumber \\
 & +& {\nb }^{i_1}{\nb }_{j_1}
 ...{\nb }^{i_{l-1}}{\nb }_{j_{l-1}}({\nb }^s {\nb }_{j_l} v) R^{i_l}_{\;\;ski} \\
 & & \nonumber \\
 &-&{\nb }^{i_1} {\nb }_{s}({\nb }^{i_2}{\nb }_{j_2}
 ...{\nb }^{i_l}{\nb }_{j_l} v )R^{s}_{\;\;j_1 ki}-...\nonumber \\
 & & \nonumber \\
 &-&{\nb }^{i_1}{\nb }_{j_1}
 ...{\nb }^{i_{l-1}}{\nb }_{j_{l-1}}({\nb }^{i_l} {\nb }_{s} v)R^{s}_{\;\;j_l ki}.\nonumber
 \end{eqnarray}
 Then we put $i_1=j_1$,...,$i_l=j_l$ in (\ref{d2}), sum up and
 cancel $l$ pairs of terms involving the Riemann tensor. Thus
 \bb\label{d3} {\nb }_i{\nb }_k ({\De }^l_g v ) ={\nb }_k{\nb }_i({\De }^l_g v
 ).\ee
 Further, we choose $p=l+1$, $q=l$, substitute
 \[T^{i i_1...i_l}_{\;\;\;\;\;\;\;\;\;\;\;\;j_1...j_l}={\nb }^i{\nb }^{i_1}{\nb }_{j_1}
 ...{\nb }^{i_l}{\nb }_{j_l} v\]
 into (\ref{d1}), and sum up over $i_1=j_1$,...,$i_l=j_l$. In this
 way we get:
 \[{\nb }_i {\nb }_k {\nb }^i ({\De }_g^l v)= {\nb }_k {\nb }_i {\nb }^i ({\De }_g^l
 v) + {\nb }^s({\De }_g^l v)R^{i}_{\;ski}.\]
 Hence
 \bb\label{d4}{\nb }_i {\nb }_k {\nb }^i ({\De }_g^l v)= {\nb }_k ({\De }_g^{l+1} v)+R^s_{\;k} {\nb }_s ({\De }_g^l
 v).\ee
 Then from (\ref{d3}) and (\ref{d4}) it follows that
 \begin{eqnarray}
 {\De }_g ({\nb }_k ({\De }_g^l v))& =&{\nb }_i {\nb }^i{\nb }_k ({\De }_g^l v) =
 {\nb }_i {\nb }_k {\nb }^i ({\De }_g^l v)\nonumber \\
 & & \nonumber \\
 & = &{\nb }_k ({\De }_g^{l+1} v)+R^s_{\;k} {\nb }_s ({\De }_g^l
 v).\nonumber \end{eqnarray}
 This completes the proof of Lemma 4.4.\end{proof}

Now let us suppose that $M$ admits a conformal Killing vector field
$h = h^i \frac{\kd }{\kd x^i}$ such that \[ {\nabla }^k {h }^s +
{\nabla }^s {h }^k = c \;g^{ks} =\frac{2}{n}\; div (h ) \;g^{ks},
\]
 where $c\neq 0$ is a constant. That is,
we suppose that $M$ admits a homothety which is not an
infinitesimal isometry of $M$. In this case we may assume that
$c=2$  and hence \bb\label{n23} div (h ) = n.\ee (Otherwise, since
$c\neq 0$, we could consider $2 h /c$ instead of $h $.)

\begin{prop} {\it   Let $u,v\in C^{2l+1}({M})$ be  given
 functions and let  $h =h^i(x) \frac{\kd }{\kd x^i}$ be 
 a $C^1({M})$ conformal Killing vector field such that $div(h)=n$.
 Then for any $l\in \mathbb{N}$ we have}

 \bb\label{n24} {\De }_g^l \eta =2l\; {\De }_g^l v + h^k{\nb }_k({\De
 }_g^l v),\ee
 \bb\label{n25} {\De }_g^l \phi =2l\; {\De }_g^l u + h^k{\nb }_k({\De
 }_g^l u).\ee
\end{prop}

\begin{proof} We shall only prove (\ref{n24}). We shall use an
 induction argument on $l$.

 1.) Let $l=1$. Differentiating two times $ \eta = (h,\nb v)=h^k
 v_k$ we obtain
 \[ {\nb }^i\eta ={\nb }^i h^k . v_k + h^k{\nb }_k {\nb }^i v, \]
 \bb\label{n26} {\De }_g\eta = {\nb }_i{\nb }^i\eta = {\De }_g h^k.v_ k +2{\nb }^ih^k .
 {\nb }_i{\nb }_k v + h^k{\nb }_i{\nb }_k{\nb }^i v.\ee

 Then (\ref{n24}) follows from (\ref{n20}), (\ref{n21}) with $\mu
 =2$, $l=1$,
 (\ref{n22}) with $l=0$ and (\ref{n26}).

 2.) Now we suppose that (\ref{n24}) holds for some
 $l\in\mathbb{N}$. We have to prove that (\ref{n24})
 holds for
 $l+1$. Differentiating two times (\ref{n24}) we have
 \begin{eqnarray}
\quad\quad{\De }_g^{l+1}\eta & =&{\nb }_i {\nb }^i ({\De }_g^l \eta) \nonumber \\ 
& & \nonumber \\
&=& {\nb }_i [\;2l\; {\nb }^i ({\De }_g^l v) + {\nb }^ih^k .
 {\nb }_k ({\De }_g^l v) + h^k{\nb }^i {\nb }_k ({\De }_g^l v)\;]\nonumber \\
 & & \nonumber \\
 &=&2l\;({\De }_g^{l+1} v)+{\De }_g h^k.{\nb }_k ({\De }_g^l v)\nonumber \\
 & & \nonumber \\ 
 &+& 2{\nb }^ih^k .{\nb }_i {\nb }_k ({\De }_g^l v) + h^k {\De }_g {\nb }_k({\De }_g^l v)\nonumber \\
 & & \nonumber \\
 & = &2l\;({\De }_g^{l+1} v)-R^k_{\;s}h^s.{\nb }_k({\De }_g^l v)
 +2{\De }_g^{l+1}v\nonumber \\
 & & \nonumber \\
 & + &h^k{\nb }_k({\De }_g^{l+1} v)+R^s_{\;k}h^k.{\nb }_s({\De }_g^l v)\nonumber \\
 & & \nonumber \\
 &=&2(l+1)\;({\De }_g^{l+1} v)+h^k{\nb }_k({\De }_g^{l+1} v).\nonumber \end{eqnarray}
 
 In the computations above we have used (\ref{n20}), (\ref{n21}) with $\mu =2$, and
 (\ref{n22}). This completes the proof. \end{proof}

Further, from (\ref{n18}), (\ref{n24}) and
 (\ref{n25}) we can express $w^i$ as
 \begin{eqnarray}\label{n27}
 w^i  = &-&\sum_{l=0}^{m-1}(2l\; {\De }_g^l v + h^k{\nb }_k({\De
 }_g^l v))\;{\nb }^i ({\De }^{2m-1-l}_g
 u) \nonumber \\
 & & \nonumber \\
 &+&\sum_{l=0}^{m-1}(2l\;{\nb }^i ({\De }^{l}_g v) +{\nb }^i h^k {\nb
 }_k ({\De }^{l}_g v) + h^k {\nb }^i{\nb }_k ({\De }^{l}_g v))
 \;({\De }^{2m-1-l}_g u)\nonumber \\
 & & \\
 &-&\sum_{l=0}^{m-1}(2l\; {\De }_g^l u + h^k{\nb }_k({\De
 }_g^l u))\;{\nb }^i ({\De }^{2m-1-l}_g
 v) \nonumber \\
 & & \nonumber \\
 &+&\sum_{l=0}^{m-1}(2l\;{\nb }^i ({\De }^{l}_g u) +{\nb }^i h^k {\nb
 }_k ({\De }^{l}_g u) + h^k {\nb }^i{\nb }_k ({\De }^{l}_g u))
 \;({\De }^{2m-1-l}_g v).\nonumber\end{eqnarray}

\begin{prop} {\it   Let $u,v\in C^{4m}(\bar{M})$ be given
 functions and let  $h =h^i(x) \frac{\kd }{\kd x^i}$ be a
 $C^1(\bar{M})$ 
  conformal Killing vector field such that $div(h)=n$.
 Then the following identity holds:}
\begin{eqnarray}\label{n28}
 R_{2m}(u,v)&=&\ds{\int_M} \{ {\De }_g^{2m} u\;(h,\nb v) +
 {\De }_g^{2m} v\;(h,\nb u)\} dV \nonumber \\
 & &\nonumber \\
 &=& (4m-n)\ds{\int_M}\lm u\lm v\; dV \\
 & & \nonumber\\
 & + &\ds{\int_{\kd M}}[ \lm u\lm v (h,\nu )-
 (w,\nu )]dS,\nonumber \end{eqnarray}
 {\it where $w=(w^i)$ is given in} (\ref{n27}).
\end{prop}

\begin{proof} We substitute into (\ref{n17}) $\lm\eta $ and
 $\lm\phi $ from (\ref{n24}) and (\ref{n25}) respectively. In this way we
 obtain:
 \[ {\De }_g^{2m} u\;(h,\nb v) +
 {\De }_g^{2m} v\;(h,\nb u)= 4m\; \lm u\lm v + h^k {\nb }_k (\lm u\lm v) - {\nb }_iw^i.\]
 Hence
 \begin{eqnarray}
  R_{2m}(u,v)&=&4m \ds{\int_M}\lm u\lm v\;dV + \ds{\int_M} h^k {\nb }_k (\lm u\lm
  v)\;dV \nonumber \\ & & \nonumber \\ &-&\ds{\int_{\kd M}}(w,\nu )dS \nonumber \\
  & & \nonumber \\
  &=&4m \ds{\int_M}\lm u\lm v\; dV-\ds{\int_M}div(h) (\lm u\lm
  v)\;dV\nonumber \\
  & & \nonumber \\
  &+ &\ds{\int_{\kd M}}[ \lm u\lm v (h,\nu )-
 (w,\nu )]dS,\nonumber \end{eqnarray}
 which implies (\ref{n28}) recalling that $div(h)=n$ (see
 (\ref{n23})).
 \end{proof}

After $2m$ integrations by parts, the first term in the
 right-hand side of (\ref{n28}) can be written in the following
 form.

\begin{prop}{\it   Let $u,v\in C^{4m}(\bar{M}),$ then the following identity holds:}
 \begin{eqnarray}\label{n29}
 \ds{\int_M}\lm u\lm v\; dV & = & \ds{\int_M} {\De }_g^{2m}
 u\;v\;dV +\sum_{l=0}^{m-1}\ds{\int_{\kd M}}{\De }_g^{m+l}u (\nb
 ({\De }_g^{m-l-1} v),\nu )dS
 \nonumber \\
 & & \\
 &-& \sum_{l=0}^{m-1}\ds{\int_{\kd M}}{\De }_g^{m-l-1}v\; (\nb
 ({\De }_g^{m+l} u),\nu )\;dS. \nonumber \end{eqnarray}
 \end{prop}

 As a consequence,  the identity (\ref{n30}) follows from (\ref{n28}) and (\ref{n29}). This completes the proof of Theorem 1.3. 

\section{A biharmonic Rellich type identity}
 
In this section we prove Theorem 1.5.
Let $\eta $ and $\phi $ be as in the preceding section. However,
 let $h$ be a conformal Killing vector field satisfying
 (\ref{n19}), where the function $\mu =2\; div(h)/n$ is not
 necessarily a constant.

 The identity (\ref{n17}) with $m=1$ reads
 \bb\label{n33} {\De }_g u\; {\De }_g \eta + {\De }_g v\;{\De }_g \phi =
 {\De }^{2}_g u.\eta +{\De }^{2}_g
 v.\phi + {\nb }_i w^i, \ee
 where
 \[w^i =-\eta {\nb }^i ({\De }_g u)+{\nb }^i\eta {\De }_g u
 -\phi {\nb }^i ({\De }_g v)+{\nb }^i\phi {\De }_g v.\]
 From (\ref{n26}), (\ref{d0}), (\ref{n21}) with $l=0$ and (\ref{d4}) with $l=0$
 it follows that
 \bb\label{n34} {\De }_g \eta =\mu {\De }_g v +h^k {\nb }_k({\De }_g
 v)+\frac{2-n}{n} {\mu }^k {\nb }_k v.\ee
 Analogously
 \bb\label{n35} {\De }_g \phi=\mu {\De }_g u +h^k {\nb }_k({\De }_g
 u)+\frac{2-n}{n} {\mu }^k {\nb }_k u.\ee
 Then, substituting (\ref{n34}) and (\ref{n35}) into (\ref{n33})
 and integrating by parts, after some work, we obtain that
 \begin{eqnarray}\label{n36}
 \ds{\int_M} \{ {\De }_g^2 u\;(h,\nb v) dV &+& {\De }_g^2 v\;(h,\nb u)\} dV =
 \frac{4-n}{2}\ds{\int_M} \mu\;{\De }_g u{\De }_g v dV \nonumber
 \\
 & & \nonumber \\
 &-&\ds{\int_M} {\De }_g\mu (u{\De }_g v +v{\De }_g u)dV \nonumber
 \\
 & & \nonumber \\
 & -& \ds{\int_M}(u(\nb\mu , \nb ({\De }_g v)) + v(\nb\mu , \nb ({\De }_g u))dV\nonumber\\
 & & \\
 &+& \ds{\int_{\kd M}}({\De }_g u{\De }_g v (h,\nu ) +
 (u{\De }_g v +v{\De }_g u)(\nb\mu ,\nu) \nonumber\\
 & & \nonumber\\
 &+& (h,\nb v) {\nb }^i({\De }_g u) {\nu }_i + (h,\nb u) {\nb }^i({\De }_g v) {\nu }_i\nonumber\\
 & & \nonumber\\
 &-& {\De }_g u ({\nb }^i h^k. v_k{\nu }_i +h^k
 {\nb }_k{\nb }^i v .{\nu }_i) \nonumber \\ & &\nonumber \\ 
 &-&{\De }_g v ({\nb }^i h^k. u_k{\nu }_i +h^k
 {\nb }_k{\nb }^i u .{\nu }_i))dS. \nonumber \end{eqnarray}

 But if $h$ satisfies (\ref{n19}), then
 \bb\label{n37}{\De }_g\mu =-\frac{1}{n-1} ({\cal L}_h R+\mu R),
 \ee
 where ${\cal L}_h $ is the Lie derivative with
 respect to the vector field $h$. See \cite{yano}. From
 (\ref{n36}) and (\ref{n37}) we get (\ref{n38}).

 We observe that if the conformal factor $\mu =2$, that is, if $div(h)=n$ (see (\ref{n24})),
 then, after two integrations by parts,
 the identity (\ref{n38}) is a particular case of (\ref{n30}) with
 $m=1$.
 
\section{A nonexistence result for a higher order semilinear Hamiltonian elliptic  system}
 
In this section, we bound ourselves to point out a very simple application of the identities proved in this paper. However, as mentioned in the introduction, the interested reader can easily realize the huge number of possible applications of these identities to other related problems in a different context.
 
Consider on $(M,g)$ the following nonlinear system of two biharmonic equations
 
 \bb\label{w1}
 \left \{ \begin{array}{lll} {\De }_g^2 u&=&\ds{\frac{\kd G}{\kd v}}, \\ \\
 {\De }_g^2 v &= & \ds{\frac{\kd G}{\kd u}},
  \end{array} \right. \ee 
  with Navier boundary conditions on $\kd M$ 
  \bb\label{w2} u=v=\De u=\De v=0.\ee
  
\begin{thm}{\it Suppose that $(M,g)$ admits a $C^1(\bar{M})$
 contravariant conformal Killing vector field $h$ such that $div(h)=n$ and $(h , \nu )>0$ on $\kd M$.
 Let $G=G(s,t)\in
 C^1({\mathbb{R}}^{2})$ satisfy the conditions}

 $(1)$\ $G(0,0)=\ds{\frac{\kd G}{\kd s}(0,0)=\frac{\kd G}{\kd t}(0,0)}=
 0,$;

 $(2)$\ {\it if } $s,t\geq 0$, {\it then}
 $\ds{\frac{\kd G}{\kd s}(s,t)\geq 0}$ {\it and }
 $\ds{\frac{\kd G}{\kd t}(s,t)}\geq 0$;

 $(3)$ {\it there exist constants $c\geq  n/(n-4)$ and $a\in
 (0,1)$ such that for any $s\in {\mathbb{R}}^{1}$ and $t\in{\mathbb{R}}^{1}$}:
 \bb\label{w3} c H(s,t) \leq a s \ds{\frac{\kd G}{\kd s}} (s,t)
 + (1-a) t \ds{\frac{\kd G}{\kd t}}(s,t).\ee

 Then there is no nontrivial classical  solution $($that is $C^4(M)\cap
 C^3(\bar{M}))$ of the Hamiltonian system (\ref{w1}) with
 Navier boundary conditions.
 \end{thm}

\begin{proof} By (\ref{n30}) with $m=2$ and $l=0$, we obtain:
 \bb\label{w4} 
 \ds{\int_M} \{ {\De }_g^{2} u\;(h,\nb v) +
 {\De }_g^{2} v\;(h,\nb u)\} dV = (4-n) \ds{\int_M} {\De }_g^{2} u\;vdV +A,
 \ee
 where
 \[ A= A(u,v) = A(v,u) =\ds{\int_{\kd M}} h^k u_k {\nb }^i ({\De }_g v) {\nu }_i dS+
  \ds{\int_{\kd M}} h^k v_k {\nb }^i ({\De }_g u) {\nu }_i dS.\]
 
On the other hand, multiplying the first equation in (\ref{w1}) by $a v$, the second - by $(1-a)u$, adding and integrating by parts, taking into account the boundary conditions (\ref{w2}), we get that
\bb\label{w5} \ds{\int_M} {\De }_g^{2} u\;vdV=\ds{\int_M} {\De }_g u{\De }_g vdV = \ds{\int_M}(a u G_u+(1-a)vG_v)dV.\ee 
Then (\ref{w4}) and (\ref{w5}) imply
  \bb\label{w6} 
 \ds{\int_M} \{ {\De }_g^{2} u\;(h,\nb v) +
 {\De }_g^{2} v\;(h,\nb u)\} dV = (4-n)\ds{\int_M}(a u G_u+(1-a)vG_v)dV +A. \ee
 
Integrating by parts in the left-hand side of (\ref{w6}) and using $div(h)=n$ we obtain that 
 \bb\label{w7} 
  n\ds{\int_M} G(u,v)dV= (n-4)\ds{\int_M}(a u G_u+(1-a)vG_v)dV -A. \ee
  
The conditions on $G$ and the maximum principle imply that 
  \[ \frac{\kd u}{\kd \nu}<0, \frac{\kd v}{\kd \nu}<0, \frac{\kd ({\De }_g u)}{\kd \nu}>0, \frac{\kd ({\De }_g v)}{\kd \nu}>0\]
  on $\kd M$. Hence
  \[0>\frac{\kd u}{\kd \nu} \frac{\kd ({\De }_g v)}{\kd \nu}=g^{ik}u_k {\nu }_i g^{js} ({\De }_g v)_s {\nu }_j=
  g^{ik}g^{js}u_j{\nu }_i ({\De }_g v)_s {\nu }_k\] 
  \[\;\;\;\;\;\;\;\;=g^{js}u_j ({\De }_g v)_s =(\nb u, \nb ({\De }_g v)), \]
  that is
  \bb\label{w8} (\nb u, \nb ({\De }_g v))<0\ee
  on $\kd M$. Above we have used the fact that on the boundary \bb\label{v1} u_k{\nu }_j= u_j{\nu }_k\ee for (\ref{w2}) 
  (see \cite{em}) and also $|\nu |^2=g^{ik}{\nu }_i {\nu }_i =1$. 
Analogously
  \bb\label{w9} (\nb v, \nb ({\De }_g u))<0\ee
on $\kd M$. Further 
 \bb\label{v2} A=  \ds{\int_{\kd M}} (\nb u, \nb ({\De }_g v)) (h,\nu ) dS +
 \ds{\int_{\kd M}} (\nb v, \nb ({\De }_g u)) (h,\nu ) dS.\ee 
Then (\ref{w8}), (\ref{w9}), (\ref{v2}) and $(h,\nu )>0$ imply $-A>0$. Then from (\ref{w7}) it follows that 
 \[  \frac{n}{n-4}\ds{\int_M} G(u,v)dV > (n-4)\ds{\int_M}(a u G_u+(1-a)vG_v)dV \]
which contradicts (\ref{w3}). This completes the proof. \end{proof} 

\subsection*{Acknowledgment}

We kindly thank the referee for the
careful reading of the paper and for pointing out how to improve its readability.

 \end{document}